\documentclass[12pt]{article}

\usepackage{amsfonts}
\usepackage{amsxtra}
\usepackage{amsthm}
\usepackage{pstricks}

\newcommand{\vect}[1]{\mathbf{#1}}
\newcommand{\abs}[1]{{\mathopen\mid}{#1}{\mathclose\mid}}
\begin{document}

\title{Teaching the Kepler laws for freshmen.}

\author{Maris van Haandel\thanks{Supported by NWO.} \\
RSG Pantarijn, Wageningen \\
marisvanhaandel@wanadoo.nl
\and 
Gert Heckman \\
IMAPP, Radboud Universiteit, Nijmegen \\
G.Heckman@math.ru.nl}

\maketitle

\begin{abstract}
We present a natural proof of Kepler's law of ellipses in the spirit of Euclidean geometry.
Moreover we discuss two existing Euclidean geometric proofs,
one by Feynman in hist Lost Lecture from 1964 and
the other by Newton in the Principia of 1687.
\end{abstract}

\section{Introduction.}

One of the highlights of classical mechanics is the mathematical derivation of the three
experimentally observed Kepler laws of planetary motion from Newton's laws
of motion and of gravitation. 
Newton published his theory of gravitation in 1687 in the Principia Mathematica \cite{Newton}.
After two short introductions, one with definitions and the other with axioms (the laws of motion),
Newton discussed the Kepler laws in the first three sections of Book 1
(in just 40 pages, without ever mentioning the name of Kepler!)

Kepler's second law (motion is planar and equal areas are swept out in equal times)
is an easy consequence of the conservation of angular momentum
$ \vect{L}=\vect{r} \times \vect{p} $,
and holds in greater generality for any central force field.
All this is explained well by Newton in Propositions 1 and 2.

Kepler's first law (planetary orbits are ellipses with the center of the force field at a focus)
however is peculiar for the attractive $1/r^2$ force field.
Using Euclidean geometry Newton derives in Proposition 11 that the 
Kepler laws can only hold for an attractive $1/r^2$ force field.
The reverse statement that an attractive $1/r^2$ force field leads to
elliptical orbits Newton concludes in Corollary 1 of Proposition 13.
Tacitly he assumes for this argument that the equation of motion
$\vect{F}=m \vect{a}$ has a unique solution for given initial position 
and initial velocity. Theorems about existence and uniqueness of solutions
of such a differential equation have only been formulated and mathematically
rigorously been proven in the 19th century. However there can be little doubt
that Newton did grasp these properties of his equation 
$\vect{F}=m \vect{a}$ \cite{Arnold}.

Somewhat later in 1710 Jakob Hermann and Johan Bernoulli gave a direct proof of Kepler's first law,
which is still the standard proof for modern text books on classical mechanics \cite{Speiser}.
One writes the position vector $\vect{r}$ in the plane of motion
in polar coordinates $r$ and $\theta$.
The trick is to transform the equation of motion $m\vect{a} = -k \vect{r}/r^3$ 
with variable the time $t$ into a second order differential equation
of the scalar function $u=1/r$ with variable the angle $\theta$.
This differential equation can be exactly solved,
and yields the equation of an ellipse in polar coordinates \cite{Goldstein1}.

Another popular proof goes by writing down the so called Laplace-Runge-Lenz vector
\[ \vect{K} = \vect{p} \times \vect{L} - km \vect{r}/r \]
with $\vect{p}=m\vect{v}$ momentum and 
$\vect{F}(\vect{r}) = -k \vect{r}/r^3$ the force field of the Kepler problem.
The LRL vector $\vect{K}$ turns out to be conserved, i.e. 
$\dot{\vect{K}}=\vect{0}$. 
This result can be derived by a direct computation as indicated in Section 2.
An alternative geometric argument is sketched in Section 3.
Working out the equation
$\vect{r} \cdot \vect{K} = rK \cos{\theta}$
gives the equation of an ellipse in polar coordinates \cite{Goldstein1}.
The geometric meaning of the LRL vector becomes clear only a posteriori 
as a vector pointing in the direction of the long axis of the ellipse.
But the start of the proof to write down the LRL vector remains a trick.
For a historical account of the LRL vector we refer to Goldstein \cite{Goldstein2}.

The purpose of this note is to present in Section 2 a proof of the Kepler laws
for which a priori the reasoning is well motivated in both physical and geometric terms.
In Section 3 we review the hodographic proof as given by
Feynman in his "Lost Lecture" \cite{Goodstein}, 
and in Section 4 we discuss Newton's proof from the Principia\cite{Newton}.
All three proofs are based on Euclidean geometry,
although we do use the language of vector calculus
in order to make it more readable for people of the 21st century.
We feel that our proof is really the simplest of the three,
and at the same time it gives more refined information
(namely the value of the long axis $2a=-k/H$ of the ellips $\mathcal{E}$).
In fact we think that our proof in Section 2 can compete equally well
both in transparency and in level of computation
with the standard proof of Jakob Hermann and Johann Bernoulli, 
making it an appropriate alternative to present in a freshman course on classical mechanics.

We like to thank Hans Duistermaat and Arnoud van Rooij for useful comments on this article.

\section{A Euclidean proof of Kepler's first law.}

We shall use inner products $\vect{u} \cdot \vect{v}$ 
and outer products $\vect{u} \times \vect{v}$ of vectors
$\vect{u}$ and $\vect{v}$ in $\mathbb{R}^3$
with the compatibility conditions
\[ \vect{u} \cdot (\vect{v} \times \vect{w}) = (\vect{u} \times \vect{v}) \cdot \vect{w} \]
\[ \vect{u} \times (\vect{v} \times \vect{w}) = (\vect{u} \cdot \vect{w}) \vect{v} - (\vect{u} \cdot \vect{v}) \vect{w} \]
and the Leibniz product rules
\[(\vect{u} \cdot \vect{v}) \spdot = \dot{\vect{u}} \cdot \vect{v} + \vect{u} \cdot \dot{\vect{v}} \]
\[(\vect{u} \times \vect{v}) \spdot = \dot{\vect{u}} \times \vect{v} + \vect{u} \times \dot{\vect{v}} \]
without further explanation.

For a central force field $\vect{F}(\vect{r}) = f(\vect{r})\vect{r}/r$  
the angular momentum vector $\vect{L}=\vect{r} \times \vect{p}$ is conserved by 
Newton's law of motion $\vect{F}=\dot{\vect{p}}$, thereby leading to Kepler's second law.
For a spherically symmetric central force field $\vect{F}(\vect{r}) = f(r)\vect{r}/r$
the energy
\[ H = p^2/2m + V(r) \; , \; V(r) = - \int f(r)\,dr \]
is conserved as well. 
These are the general initial remarks.

From now on consider the Kepler problem $f(r)=-k/r^2$ en $V(r)=-k/r$ with
$k>0$ a coupling constant. For fixed energy $H<0$ the motion is bounded
inside a sphere with center $\vect{0}$ and radius $-k/H$.
Consider the following picture of the plane perpendicular to $\vect{L}$.

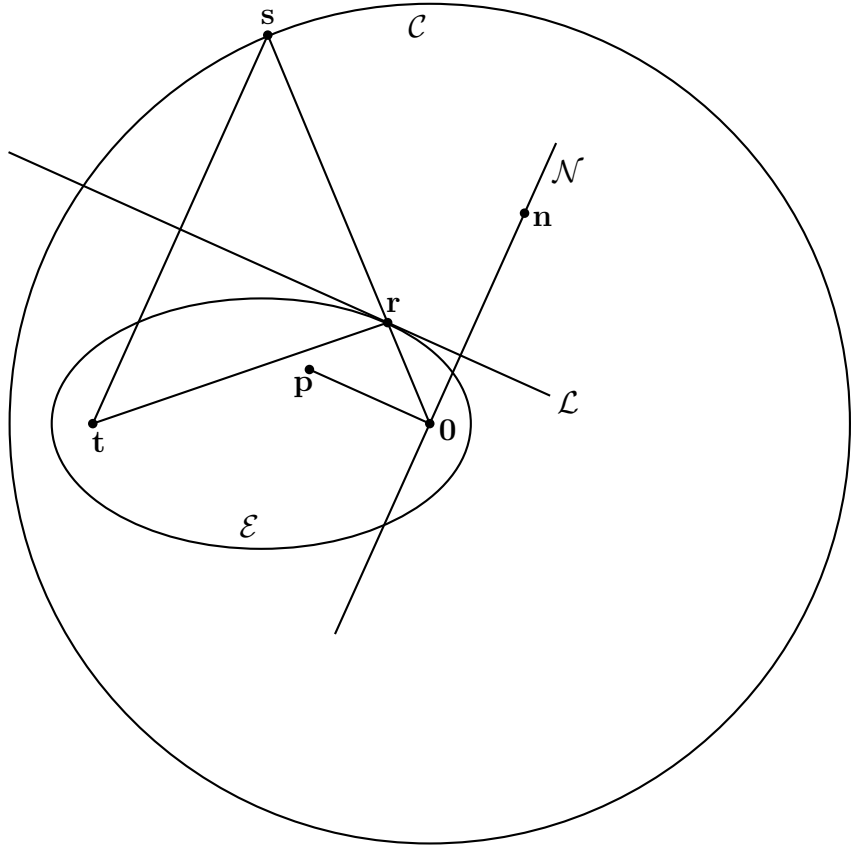
\begin{center}
\psset{unit=0.8mm}

\begin{pspicture}*(-75,-75)(75,75)

\pscircle(0,0){70}
\psellipse(-28,0)(35,21)
\psdot(0,0)
\psdot(-56,0)
\psdot(-20,9)
\psdot(-7,16.8)
\psdot(15.75,35)
\psdot(-26.923,64.615)
\psline(-70,45.15)(20,4.65)
\psline(-15.75,-35)(21,46.666)
\psline(-26.923,64.615)(-56,0)
\psline(-7,16.8)(-56,0)
\psline(0,0)(-20,9)
\psline(0,0)(-26.923,64.615)
\rput(3,-1){$\vect{0}$}
\rput(18.75,34){$\vect{n}$}
\rput(-55,-3){$\vect{t}$}
\rput(-21,6){$\vect{p}$}
\rput(-6,19.8){$\vect{r}$}
\rput(-26.923,67.615){$\vect{s}$}
\rput(-2,66){$\mathcal{C}$}
\rput(-30,-17){$\mathcal{E}$}
\rput(23,3.65){$\mathcal{L}$}
\rput(23,42){$\mathcal{N}$}

\end{pspicture}
\end{center}

The circle $\mathcal{C}$ with center $\vect{0}$ and radius $-k/H$ is the boundary of
a disc where motion with energy $H<0$ takes place.
Let $\vect{s}=-k \vect{r}/rH$ be the projection of $\vect{r}$ from the center
$\vect{0}$ on this circle $\mathcal{C}$. 
The line $\mathcal{L}$ through $\vect{r}$ with direction vector $\vect{p}$ is the tangent line
of the orbit $\mathcal{E}$ at position $\vect{r}$ with velocity $\vect{v}$.
Let $\vect{t}$ be the orthogonal reflection of the point $\vect{s}$ in the line $\mathcal{L}$.

\newtheorem*{tconservedtheorem}{Theorem}

\begin{tconservedtheorem}
The point $\vect{t}$ equals $\vect{K}/mH$ and therefore is conserved.
\end{tconservedtheorem}

\begin{proof}
The line $\mathcal{N}$ spanned by $\vect{n}=\vect{p} \times \vect{L}$
is perpendicular to $\mathcal{L}$.
The point $\vect{t}$ is obtained from $\vect{s}$ by subtracting twice
the orthogonal projection of $\vect{s} - \vect{r}$ on the line $\mathcal{N}$, 
and therefore
\[ \vect{t}=\vect{s}-2((\vect{s} - \vect{r})\cdot\vect{n})\vect{n}/n^2. \]
Now
\[ \vect{s}=-k \vect{r}/rH \]
\[(\vect{s}-\vect{r})\cdot\vect{n}=-(H+k/r)\vect{r}\cdot(\vect{p}\times\vect{L})/H=-(H+k/r)L^2/H\]
\[ n^2=p^2 L^2=2m(H+k/r)L^2 \]
and therefore
\[ \vect{t} = -k\vect{r}/rH + \vect{n}/mH = \vect{K}/mH \]
with $\vect{K} = \vect{p} \times \vect{L} - km \vect{r}/r$ the LRL vector.
The final step $\dot{\vect{K}}=\vect{0}$ is derived by a straightforward computation
using the compatibility relations and the Leibniz product rules
for inner and outer products of vectors in $\mathbb{R}^3$.
\end{proof}

\newtheorem*{ellipsecorollary}{Corollary}

\begin{ellipsecorollary} 
The orbit $\mathcal{E}$ is an ellipse with foci
$\vect{0}$ and $\vect{t}$, and long axis equal to $2a=-k/H$.
\end{ellipsecorollary}

\begin{proof}
Indeed we have
\[ \abs{\vect{t}-\vect{r}} + \abs{\vect{r}-\vect{0}} =
   \abs{\vect{s}-\vect{r}} + \abs{\vect{r}-\vect{0}} =
   \abs{\vect{s}-\vect{0}} = -k/H. \]
Hence $\mathcal{E}$ is an ellipse with foci 
$\vect{0}$ and $\vect{t}$, and long axis $2a=-k/H$.
\end{proof}

The above proof has two advantages over the earlier mentioned proofs of Kepler's first law.
The conserved vector $\vect{t} = \vect{K}/mH$ is a priori well motivated in geometric terms. 
Moreover we use the gardener's definition of an ellipse.
The equation of an ellipse in polar coordinates is unknown to most freshmen,
and so additional explanation would be needed for that.
Another advantage of our proof is that the solution of the equation of motion
is achieved by just finding enough constants of motion.
The proofs by Feynman and Newton in the next sections on the contrary rely
at a crucial point on the existence and uniqueness theorem of differential equations.

We proceed to derive Kepler's third law along standard lines \cite{Goldstein1}.
The ellipse $\mathcal{E}$ has numerical parameters\footnote{The long axis equals $2a$,
the short axis $2b$, and $a^2=b^2+c^2$.} $a,b,c>0$ given by
$2a=-k/H$, $4c^2=K^2/m^2H^2=(2mHL^2+m^2k^2)/m^2H^2$ and $a^2=b^2+c^2$.
The area of the region bounded by $\mathcal{E}$ equals
\[ \pi ab = LT/2m \]
with $T$ the period of the orbit. 
Indeed L/2m is the area of the sector swept out by the position vector $\vect{r}$ per unit time.
A straightforward calculation gives
\[ T^2/a^3 = 4 \pi^2 m/k. \]
Newton's law of gravitation states that the coupling constant $k$ is proportional to the product of the 
mass $m$ of the planet and the mass $M$ of the sun, 
and we find Kepler's third (harmonic) law
stating that $T^2/a^3$ is the same for all planets\footnote{The mass m we have used
so far is in fact equal to the reduced mass $mM/(m+M)$,
and this almost equals $m$ if $m \ll M$. 
The coupling constant k is according to Newton equal to $GmM$ with $G$ the universal gravitational constant.
We therefore see that Kepler's third law holds only approximately for $m \ll M$.}.

\section{Feynman's proof of Kepler's first law.}

In this section we discuss a different geometric proof of Kepler's first law
based on the hodograph $\mathcal{H}$. By definition $\mathcal{H}$ is the curve
traced out by the velocity vector $\vect{v}$ in the Kepler problem.
This proof goes back to M\"obius in 1843 and Hamilton in 1845 \cite{Derbes}
and has been forgotten and rediscovered several times,
among others by Feynman in 1964 in his "Lost Lecture" \cite{Goodstein}.

Let us assume (as in the picture of the previous section)
that $iv\vect{n}/n=\vect{v}$ with $i$ the counterclockwise
rotation around $\vect{0}$ over $\pi/2$. 
So the orbit $\mathcal{E}$ is assumed to be traversed
counterclockwise around the origin $\vect{0}$.

\newtheorem*{hodotheorem}{Theorem}

\begin{hodotheorem}
The hodograph $\mathcal{H}$ is a circle with center
$\vect{c}=i\vect{K}/mL$ and radius $k/L$
\end{hodotheorem}

\begin{proof}
We shall indicate two proofs of this theorem. 
The first proof is analytic in nature, and uses
conservation of the LRL vector $\vect{K}$ by rewriting
\[ \vect{K} = \vect{p} \times \vect{L} - km \vect{r}/r = mvL \vect{n}/n -km \vect{r}/r \]
as
\[ v \vect{n}/n = \vect{K}/mL + k \vect{r}/rL \;, \]
or equivalently
\[ \vect{v} = i \vect{K}/mL + ik \vect{r}/rL \;. \]
Hence the theorem follows from $\dot{\vect{K}}=\vect{0}$.

There is a different geometric proof of the theorem, which as a
corollary gives the conservation of the LRL vector $\vect{K}$.
The key point is to reparametrize the velocity vector $\vect{v}$
from time $t$ to angle $\theta$ of the position vector $\vect{r}$.
It turns out that the vector $\vect{v}(\theta)$ is traversing the
hodograph $\mathcal{H}$ with constant speed $k/L$.
Indeed we have from Newton's laws
\[ m \frac{d \vect{v}}{d \theta} = m \vect{a} \frac{dt}{d \theta} = - \frac{k \vect{r}}{r^3} \frac{dt}{d \theta} \]
and Kepler's second law gives
\[ r^2 d \theta/2 = L dt/2m \;. \]
Combining these identities yields
\[ \frac{d \vect{v}}{d \theta} = -k \vect{r}/rL \;, \]
so indeed $\vect{v}(\theta)$ travels along $\mathcal{H}$ with
constant speed $k/L$.
Since $\vect{r} = r e^{i \theta}$ a direct integration yields
\[ \vect{v}(\theta) = \vect{c} + ik \vect{r}/rL \; , \; \dot{\vect{c}}=\vect{0} \]
and the hodograph becomes a circle with center $\vect{c}$ and radius $k/L$.
Comparison with the last formula in the first proof gives
\[ \vect{c} = i \vect{K}/mL \]
and $\dot{\vect{K}}=\vect{0}$ comes out as a corollary.
\end{proof}

All in all, the circular nature of the hodograph $\mathcal{H}$
is more or less equivalent to the conservation of the LRL vector $\vect{K}$.

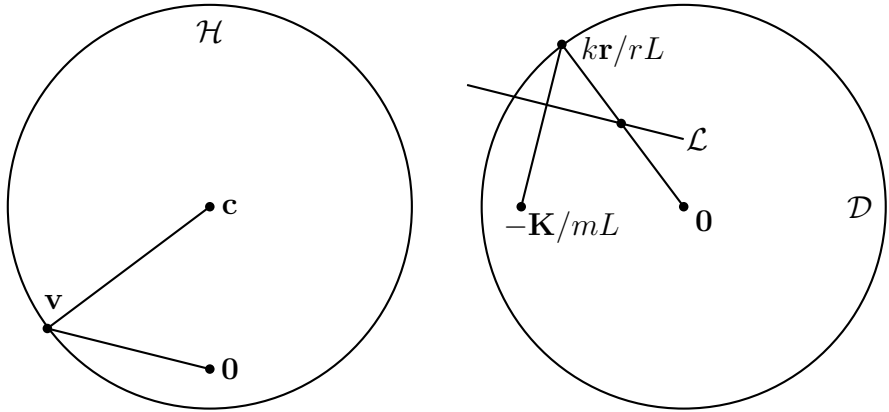
\begin{center}
\psset{unit=0.9mm}

\begin{pspicture}*(-75,-35)(75,35)

\pscircle(-35,0){30}
\pscircle(35,0){30}
\psdot(-35,0)
\psdot(35,0)
\psdot(-35,-24)
\psdot(11,0)
\psdot(-59,-18)
\psdot(17,24)
\psdot(25.769,12.308)%
\psline(-35,-24)(-59,-18)
\psline(11,0)(17,24)
\psline(17,24)(35,0)
\psline(-59,-18)(-35,0)
\psline(3,18)(35,10)
\rput(-32,-24){$\vect{0}$}
\rput(38,-2){$\vect{0}$}
\rput(-32,0){$\vect{c}$}
\rput(-58,-14){$\vect{v}$}
\rput(17,-3){$-\vect{K}/mL$}
\rput(26,23){$k \vect{r}/rL$}
\rput(-35,26){$\mathcal{H}$}
\rput(61,0){$\mathcal{D}$}
\rput(37,10){$\mathcal{L}$}

\end{pspicture}
\end{center}

Now turn the hodograph $\mathcal{H}$ clockwise around $\vect{0}$ over $\pi /2$
and translate over $i \vect{c}=- \vect{K}/mL$. 
This gives a circle $\mathcal{D}$ with center $\vect{0}$ and radius $k/L$. Since
\[ k \vect{r}/rL + \vect{K}/mL = v \vect{n}/n = -i \vect{v}\]
the orbit $\mathcal{E}$ intersects the line through $\vect{0}$ and $k \vect{r}/rL$
in a point with tangent line $\mathcal{L}$ perpendicular to the line through
$k \vect{r}/rL$ and $- \vect{K}/mL$.
For example the ellipse $\mathcal{F}$ with foci $\vect{0}$ and $- \vect{K}/mL$
and long axis equal to $k/L$ has this property,
but any scalar multiple $\lambda \mathcal{F}$ with $\lambda > 0$ equally does so.
Since curves with the above property are uniquely charcterized once an
initial point on the curve is chosen we conclude that $\mathcal{E} = \lambda \mathcal{F}$
for some $\lambda > 0$. This proves Kepler's first law.
A comparison with the picture of the previous section shows that 
$\mathcal{E} = \lambda \mathcal{F}$ with $\lambda=-L/H$. 
Indeed $\mathcal{E}$ has foci $\vect{0}$ and 
$- \lambda \vect{K}/mL = \vect{K}/mH = \vect{t}$ and
the long axis is equal to $\lambda k/L = -k/H = 2a$.

It is not clear to us if Feynman was aware of his use of the existence and uniqueness
theorem for differential equations. On page 164 of \cite{Goodstein} he writes:
"Therefore, the solution to the problem is an ellipse - or the other way around,
really, is what I proved: that the ellipse is a possible solution to the problem.
And it is this solution. So the orbits are ellipses."

Apparently Feynman had trouble following Newton's proof of Kepler's first law. 
On page 111 of \cite{Goodstein} the authors write:
"In Feynman's lecture, this is the point at which he finds himself unable
to follow Newton's line of argument any further, and so sets out to invent one of his own".

\section{Newton's proof of Kepler's first law.}

In this section we discuss the original proof by Newton
of Kepler's first law as given in \cite{Newton}. 
The proof starts with a nice general result.

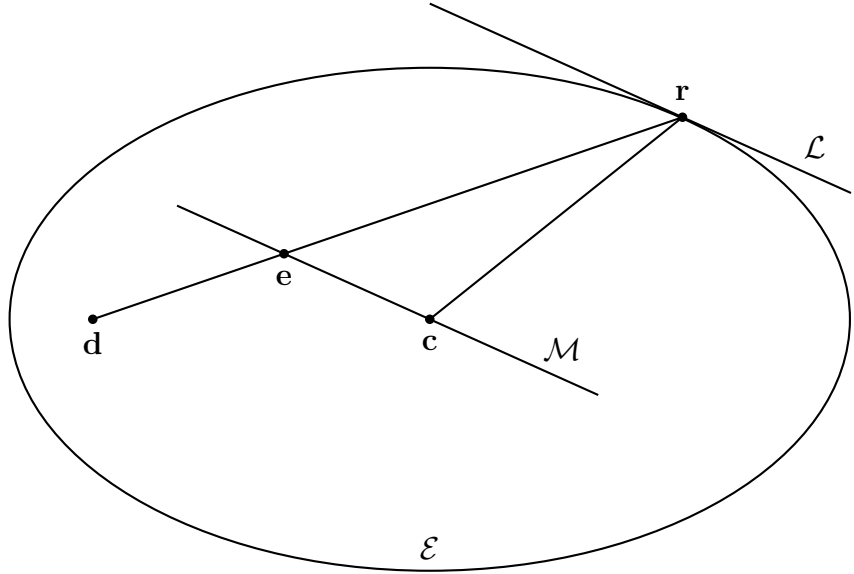
\begin{center}
\psset{unit=0.8mm}

\begin{pspicture}*(-75,-50)(75,60)

\psellipse(0,0)(70,42)
\psdot(0,0)
\psdot(-56,0)
\psdot(42,33.6)
\psdot(-24.216,10.897)
\psline(0,52.5)(70,21)
\psline(-42,18.9)(28,-12.6)
\psline(42,33.6)(0,0)
\psline(42,33.6)(-56,0)
\rput(0,-4){$\vect{c}$}
\rput(-56,-4){$\vect{d}$}
\rput(42,37.6){$\vect{r}$}
\rput(-24.216,6.897){$\vect{e}$}
\rput(0,-38){$\mathcal{E}$}
\rput(64,28.6){$\mathcal{L}$}
\rput(22,-5){$\mathcal{M}$}

\end{pspicture}
\end{center}

\newtheorem*{newtontheorem}{Theorem}

\begin{newtontheorem}
Let $\mathcal{E}$ be a smooth closed curve bounding a convex region 
containing two points $\vect{c}$ and $\vect{d}$. 
Let $\vect{r}(t)$ traverse the curve $\mathcal{E}$ counterclockwise in time $t$, 
such that the areal speed with respect to the point $\vect{c}$ is constant.
Likewise let $\vect{r}(s)$ traverse the curve $\mathcal{E}$ counterclockwise in time $s$, 
such that the areal speed with respect to the point $\vect{d}$ is equal to the same constant.

Let $\mathcal{L}$ be the tangent line to $\mathcal{E}$ at the point $\vect{r}$,
and let $\vect{e}$ be the intersection point of the line $\mathcal{M}$,
which is parallel to $\mathcal{L}$ through the point $\vect{c}$, 
and the line through the points $\vect{r}$ and $\vect{d}$.
Then the ratio of the two accelerations is given by
\[ \abs{\frac{d^2 \vect{r}}{ds^2}} : \abs{\frac{d^2 \vect{r}}{dt^2}} = 
{\abs{\vect{r}-\vect{e}}}^3 : {\abs{\vect{r}-\vect{c}}}\cdot{\abs{\vect{r}-\vect{d}}}^2 .\]
\end{newtontheorem}

\begin{proof}
Using the chain rule we get
\[ \frac{d \vect{r}}{ds} = \frac{d \vect{r}}{dt}\cdot\frac{dt}{ds} \]
\[ \frac{d^2 \vect{r}}{ds^2} = 
\frac{d^2 \vect{r}}{dt^2}\cdot(\frac{dt}{ds})^2 + \frac{d \vect{r}}{dt}\cdot\frac{d^2t}{ds^2} \; .\]
Because 
$d^2 \vect{r}/dt^2$ is proportional to $\vect{c}-\vect{r}$
and likewise
$d^2 \vect{r}/ds^2$ is proportional to $\vect{d}-\vect{r}$
we see that
\[ \abs{\frac{d^2 \vect{r}}{ds^2}} : \abs{\frac{d^2 \vect{r}}{dt^2}} =
(\frac{dt}{ds})^2 \cdot \abs{\frac{d^2 \vect{r}}{dt^2} + \frac{d \vect{r}}{dt}\cdot\frac{d^2t}{ds^2}/(\frac{dt}{ds})^2} : 
\abs{\frac{d^2 \vect{r}}{dt^2}}= (\frac{dt}{ds})^2 \cdot \abs{\vect{r}-\vect{e}} : \abs{\vect{r}-\vect{c}} \; .\]
Since the curve $\mathcal{E}$ is traversed with equal areal speed 
relative to the two points $\vect{c}$ and $\vect{d}$ we get
\[ \abs{\frac{d\vect{r}}{dt}} \cdot \abs{\vect{r}-\vect{e}} = \abs{\frac{d\vect{r}}{ds}}\cdot\abs{\vect{r}-\vect{d}} \]
and therefore also
\[ \frac{dt}{ds} = \abs{\vect{r}-\vect{e}} : \abs{\vect{r}-\vect{d}} \; . \]
In turn this implies that
\[ \abs{\frac{d^2 \vect{r}}{ds^2}} : \abs{\frac{d^2 \vect{r}}{dt^2}} =
(\frac{dt}{ds})^2 \cdot \abs{\vect{r}-\vect{e}} : \abs{\vect{r}-\vect{c}} =
{\abs{\vect{r}-\vect{e}}}^3 : {\abs{\vect{r}-\vect{c}}}\cdot{\abs{\vect{r}-\vect{d}}}^2 \]
which proves the theorem.
\end{proof}

We shall apply this theorem in case $\mathcal{E}$ is an ellipse 
with center $\vect{c}$ and focus $\vect{d}$.
Assume that $\vect{r(t)}$ traverses the ellipse $\mathcal{E}$
in harmonic motion, say
\[ \frac{d^2 \vect{r}}{dt^2} = \vect{c} - \vect{r} \; . \]

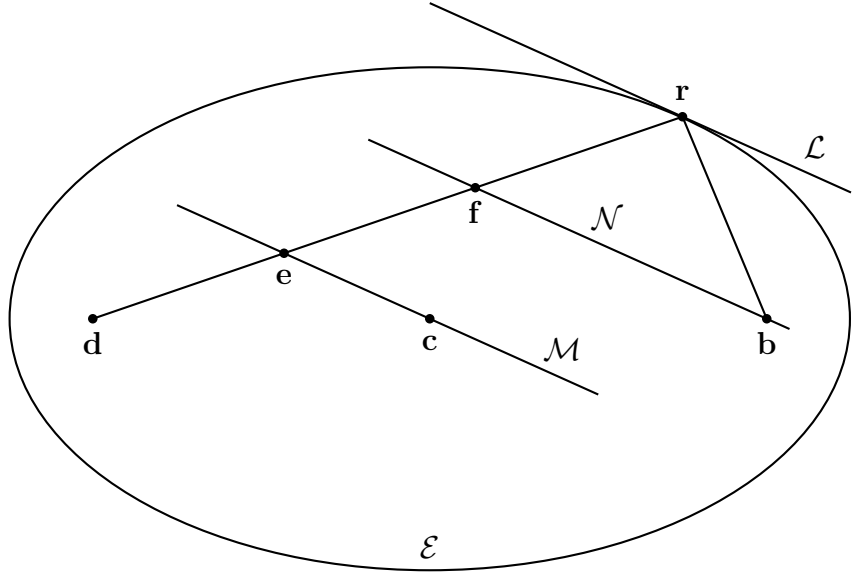
\begin{center}
\psset{unit=0.8mm}

\begin{pspicture}*(-75,-50)(75,60)

\psellipse(0,0)(70,42)
\psdot(0,0)
\psdot(-56,0)
\psdot(42,33.6)
\psdot(-24.216,10.897)
\psdot(56,0)
\psdot(7.568,21.794)
\psline(0,52.5)(70,21)
\psline(-42,18.9)(28,-12.6)
\psline(-10.216,29.797)(59.784,-1.703)
\psline(42,33.6)(56,0)
\psline(42,33.6)(-56,0)
\rput(0,-4){$\vect{c}$}
\rput(-56,-4){$\vect{d}$}
\rput(42,37.6){$\vect{r}$}
\rput(-24.216,6.897){$\vect{e}$}
\rput(56,-4){$\vect{b}$}
\rput(7.568,17.794){$\vect{f}$}
\rput(0,-38){$\mathcal{E}$}
\rput(64,28.6){$\mathcal{L}$}
\rput(22,-5){$\mathcal{M}$}
\rput(29.568,16.794){$\mathcal{N}$}

\end{pspicture}
\end{center}

Let $\vect{b}$ be the other focus of $\mathcal{E}$, 
and let $\vect{f}$ be the intersection point of the line $\mathcal{N}$,
passing through $\vect{b}$ and parallel to $\mathcal{L}$,
with the line through the points $\vect{d}$ and $\vect{r}$.
Then we find
\[ \abs{\vect{d} - \vect{e}} = \abs{\vect{e} - \vect{f}} \;, \; \abs{\vect{f} - \vect{r}} = \abs{\vect{b} - \vect{r}} \]
which in turn implies that $\abs{\vect{e} - \vect{r}}$ is equal to
the half long axis $a$ of the ellipse $\mathcal{E}$.
We conclude from the formula in the above theorem that the motion in time $s$
along an ellipse with constant areal speed with respect to a focus
is only possible in an attractive inverse square force field.
The converse statement that an inverse square force field (for negative 
energy $H$) indeed yields ellipses as orbits follows from existence and
uniqueness theorems for solutions of Newton's equation $\vect{F}=m\vect{a}$ and the above reasoning.
This is Newton's line of argument for proving Kepler's first law.

\section{Conclusion.}

There exist other proofs of Kepler's law of ellipses from a higher view point.
One such proof by Arnold uses complex analysis, and is somewhat reminiscent
to Newton's proof in Section 4 by comparing harmonic motion with motion
under an $1/r^2$ force field \cite{Arnold}. Another proof by Moser is also
very elegant, and uses the language of symplectic geometry and canonical transformations
\cite{Guillemin},\cite{Moser}. 
However our goal here has been to give a proof which is as basic as possible,
and at the same time is well motivated in terms of Euclidean geometry.

It is hard to exaggerate the importance of the role of the Principia Mathematica
in the history of science. The year 1687 marks the birth of both modern mathematical analysis
and modern theoretical physics. As such the derivation of the Kepler laws from Newton's law
of motion and law of universal gravitation is a rewarding subject to teach to freshmen students.
In fact we were motivated for our work, because we plan to teach this material to high school
students in their final grade. Of course, the high school students first need to get
acquainted with the basics of vector geometry and vector calculus. 
But once this is achieved there is nothing in the way of understanding our proof
of Kepler's law of ellipses.

For physics or mathematics freshmen students in the university who are already familiar
with vector calculus our proof in Section 2 is fairly short and geometrically well motivated.
In our opinion of all proofs it qualifies best to be discussed in a freshman course on mechanics.

\end{document}